\documentclass[11pt]{article}

\newfont{\bb}{msbm10 at 11pt}
\def\r{\hbox{\bb R}}

\newenvironment{proof}{\trivlist
\item[\hskip\labelsep{\it Proof}\,:]}{\hfill{$\Box$}\endtrivlist}
\newtheorem{theorem}{Theorem}[section]
\newtheorem{corollary}[theorem]{Corollary}
\newtheorem{lemma}[theorem]{Lemma}

\newtheorem{claim}{Claim}

\usepackage{graphicx}
\usepackage{latexsym}
\usepackage{amssymb}
 
\begin{document}

\title{Rotational linear Weingarten surfaces\\of hyperbolic type }
\author{Rafael L\'opez\footnote{Partially
supported by MEC-FEDER 
 grant no. MTM2004-00109.}\\
Departmento de Geometr\'{\i}a y Topolog\'{\i}a\\
Universidad de Granada\\
18071 Granada (Spain)\\
e-mail:{\tt rcamino@ugr.es}}
\date{}
\maketitle


\begin{abstract}
 A linear Weingarten surface in Euclidean space 
$\r^3$ is a surface 
whose mean curvature $H$ and  
Gaussian curvature $K$ satisfy a relation of the form  $aH+bK=c$, where 
$a,b,c\in\r$. Such a surface is said to be hyperbolic when  $a^2+4bc<0$. 
In this paper we classify all rotational linear Weingarten surfaces of 
hyperbolic type. As a consequence, we obtain  a family  of complete 
hyperbolic linear Weingarten surfaces in $\r^3$ that consists into  periodic 
surfaces with self-intersections.
\end{abstract}

\section{Introduction}\label{intro}

A surface $S$ in 3-dimensional Euclidean space $\r^3$ is called a {\it 
Weingarten surface} if there is some relation between its
two principal curvatures $\kappa_1$ and $\kappa_2$, 
that is, if there is a smooth function $W$ of two variables such that
$W(\kappa_1,\kappa_2)=0$. In particular, if $K$ and $H$  denote respectively  the Gauss and the mean curvature of  $S$, the identity $W(\kappa_1,\kappa_2)=0$ implies a relation $U(K,H)=0$.   The very Weingarten 
introduced this kind of surfaces in the context of the problem of finding all surfaces isometric to a given surface of revolution \cite{we1},\cite{we2}. 
In this paper we study Weingarten surfaces that satisfy the simplest case for $U$, that is, that $U$ is of linear type:
$$a\ H+b\ K=c,$$
where $a,b,c\in\r$.  We say  that $S$ is a {\it linear  Weingarten surface} and we abbreviate it by LW-surface. First examples of LW-surfaces are the surfaces with 
constant mean curvature   ($b=0$) and the surfaces with 
 constant Gauss curvature ($a=0$).  
Although these two kinds of surfaces have been extensively studied in the literature, 
the  classification of LW-surfaces in the general case  is almost completely open today. 
Along the history,
 they  have been of interest for geometers, mainly when the surface is closed: 
\cite{ch1}, \cite{ch2}, \cite{hw}, \cite{ho}, \cite{ho2}, \cite{ks}, \cite{vo}. 

The behaviour of a LW-surface and its qualitative properties strongly depends on  the 
sign of the discriminant $\Delta:=a^2+4bc$. 
Such a surface is called to be hyperbolic (resp. elliptic) when 
$\Delta<0$ (resp. $\Delta>0$). The relation $\Delta=0$ characterizes the tubular surfaces.
Examples of elliptic surfaces are the surfaces with constant mean curvature and 
the surfaces with positive constant Gaussian curvature. Since elliptic LW-surfaces have similar properties 
as these  two kinds of surfaces, they have been attracted for 
 a number of authors. For example, elliptic LW-surfaces satisfy a  
maximum principle and this enables to use Alexandrov reflection technique in its study.
See the recent bibliography \cite{be}, \cite{cft}, \cite{gmm}, \cite{rs}. 

The aim of this paper is the study of the  LW-surfaces of hyperbolic type. 
Examples of hyperbolic LW-surfaces are the surfaces with 
negative constant Gaussian curvature ($a=0, bc<0$).    One expects then 
to find in the hyperbolic LW-surfaces similar 
properties as the ones of the surfaces with negative constant Gaussian curvature. However, and as we shall see, 
our family is more extensive and richer even in the rotational case;  for example, we will obtain hyperbolic rotational 
LW-surfaces with positive Gaussian curvature (Theorem \ref{t12}). Notice that in a hyperbolic 
Weingarten surface, the condition $\Delta<0$ implies that do not exist umbilics in the surface. 
As $a^2+4bc<0$,  it 
follows that $c\not=0$. Without loss of generality, throughout this work we shall assume  
that $c=1$ and the linear Weingarten relation is now 
\begin{equation}\label{w1}
a\ H+b\ K=1.
\end{equation}
Among all hyperbolic LW-surfaces,  
 the class of the surfaces  of revolution are particularly interesting
 because  in such case, 
Equation (\ref{w1}) leads to an ordinary differential equation. Its study is then simplified to find the profile curve that defines the surface. In this paper, we classify all of rotational hyperbolic LW-surfaces. 
 We  summarize our classification as follows
 (see Theorems \ref{t12}, \ref{t2}, \ref{t32} and \ref{t42}):

\begin{quote}{\it Let $a$ and $b$ be two real numbers under the condition 
 $a^2+4b<0$. Then the rotational  linear Weingarten surfaces satisfying 
$aH+bK=1$ can parametrized by one parameter $z_0$, namely, $S(a,b;z_0)$, with $a,z_0>0$ and 
$z_0\not=-2b/a$ such that: (i) if $0<z_0<a/2$, the surface is not complete and with
positive Gaussian curvature; (ii) if $z_0=a/2$, the surface is a  right cylinder; (iii) if $a/2<z_0<-2b/a$, the surface  is not complete and with negative Gaussian curvature ; (iv) if 
$z_0>-2b/a$, the surface is complete and periodic.}
\end{quote}

\begin{figure}[htbp]\begin{center}
\includegraphics[width=9cm]{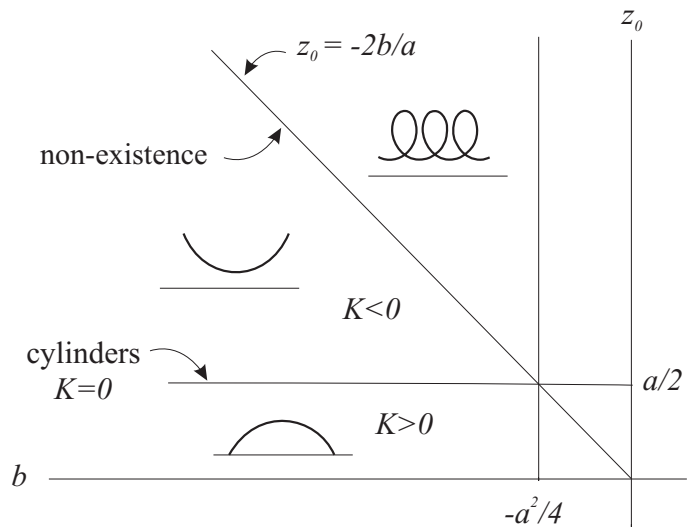}\end{center}
\caption{Classification of rotational hyperbolic LW-surfaces. We fix the value of the constant $a$ in Equation 
(\ref{w1}) with $a>0$. The diagram represents the $(b,z_0)$-plane and  let us consider $z_0$ as function of the variable $b$. When $z_0=-2b/a$, do not exist surfaces. In the other case, 
we have four families of surfaces according  the value $z_0$ goes from $0$ to $\infty$ and 
depending of its position with respect to  the values $a/2$ and $-2b/a$.}\label{fig0}
\end{figure}

In Figure \ref{fig0} we present a scheme of the classification of rotational hyperbolic 
LW-surfaces realized in this paper.  In order to get this classification, 
we shall describe the symmetries and qualitative properties of these surfaces.
Among the rotational hyperbolic LW-surfaces obtained in the above result, a special class of such surfaces is the family of surfaces $S(a,b;z_0)$ with $z_0>-2b/a$, which have remarkable properties. For this reason, we separate and emphasize  the statement  as follows  (Corollary \ref{co}):

\begin{quote} {\it There exists a one-parameter family of 
rotational hyperbolic linear Weingarten 
surfaces that 
are complete and with self-intersections in $\r^3$. Moreover, these surfaces are periodic.}
\end{quote}

This contrasts with Hilbert's theorem  that there do not exist  complete surfaces with constant negative Gaussian 
curvature immersed in $\r^3$. See also \cite{cft} for other examples. 

This paper is organized as follows. Section \ref{s0} introduces notation and terminology used throughout this 
paper. Sections \ref{s1} to \ref{s4} sucessively describe   the properties of rotational hyperbolic LW-surfaces according the value that takes the parameter $z_0$.  For each 
case depending on $z_0$, we present pictures of the profile curves.

\section{Preliminaries}\label{s0}

Let $\r^3$ be the three-dimensional space with usual coordinates $(x,y,z)$. 
Let $\alpha:I\rightarrow\r^3$ be a planar curve  in the $(x,z)$-plane with 
coordinate functions $\alpha(s)=(x(s),0,z(s))$ and $z(s)>0$. Assume that $s$ is the 
arclength along $\alpha$.  Consider $\theta=\theta(s)$ the angle 
function that makes the velocity $\alpha'(s)$ at $s$  with the 
$x$-axis, that is,  $\alpha'(s)=(\cos\theta(s),0,\sin\theta(s))$. The 
curvature of the planar curve $\alpha$ is given by $\theta'$. 
Consider $S$ the surface of revolution obtained by rotating $\alpha$ with respect to the 
$x$-axis, that is, $S$ parametrizes as $X(s,\phi)=(x(s), z(s) \cos\phi, 
z(s)\sin\phi)$.  
The principal curvatures of $S$ are given by 
$$\kappa_1(s,\phi)=\frac{\cos\theta(s)}{z(s)},\hspace*{1cm}
\kappa_2(s,\phi)=-\theta'(s),$$
and the mean curvature and the Gaussian curvature of $S$ write, respectively,  as 
\begin{equation}\label{hk}
H(s,\phi)=\frac{\cos\theta(s)-z(s)\theta'(s)}{2 z(s)},\hspace*{1cm}
K(s,\phi)=-\frac{\cos\theta(s)\theta'(s)}{z(s)}.
\end{equation}
The Weingarten relation (\ref{w1}) converts into
\begin{equation}\label{w2}
a\ \frac{\cos\theta(s)-z(s)\theta'(s)}{2 z(s)}-b\ \frac{\cos\theta(s)\theta'(s)}{z(s)}=1.
\end{equation}
Moreover, throughout this work, we
 discard that the surface has negative constant Gaussian curvature. Thus, we 
assume $a\not=0$. 

The study of rotational hyperbolic LW-surfaces  reduces to the knowledge of the solutions of Equation (\ref{w2}) for given initial data. 
Consider  that $\alpha(0)=(0,0,z_0)$, $z_0>0$ and $\alpha'(0)=(1,0,0)$.
Then the  curve $\alpha$ is governed by the differential equations system 
\begin{equation}\label{eq1}
 \left\{\begin{array}{lll}
 x'(s)&=& \displaystyle \cos\theta(s)\\
 z'(s)&=&\displaystyle \sin\theta(s)\\
 \theta'(s)&=&\displaystyle \frac{a\cos\theta(s)-2z(s)}{a z(s)+2b \cos\theta(s)}
\end{array}
\right.
\end{equation}
with initial conditions
\begin{equation}\label{eq2}
x(0)=0,\hspace*{.5cm}z(0)=z_0,\hspace*{.5cm}\theta(0)=0.
\end{equation}

If it is necessary, we shall denote by $\alpha(s;z_0)$ the solution obtained in 
(\ref{eq1})-(\ref{eq2}) to 
emphasize the dependence on the parameter $z_0$.
We have assumed that $z(s)>0$ in the parametrization of $\alpha$. However, if we change the 
generating curve $\alpha$ by $\bar{\alpha}(s)=(x(s),0,-z(s))$, the surface obtained by rotating the new curve $\bar{\alpha}$ with respect to the $x$-axis is the same 
that the $\alpha$ one. But now, the curve $\bar{\alpha}$ is the solution of (\ref{eq1})-(\ref{eq2}) changing the parameters 
$(a,b,z_0)$ by $(-a,b,-z_0)$. Hence that we can choose the parameters $a$ and $z_0$ to 
have the same sign. For convenience, we shall assume that both $a$ and $z_0$ are positive 
numbers.

\begin{theorem} A first integral of the differential equations system 
(\ref{eq1})-(\ref{eq2}) is given by 
\begin{equation}\label{first}
z(s)^2-a z(s)\cos\theta(s)-b\cos^2\theta(s)-(z_0^2-az_0-b)=0.
\end{equation}
\end{theorem}

\begin{proof}
Multiplying both sides of (\ref{w2}) by $\sin\theta$, we obtain a first integral of 
type
$$z(s)^2-a z(s)\cos\theta(s)+b\sin^2\theta(s)+\lambda=0$$
for some $\lambda\in\r$. At $s=0$, we conclude $\lambda = az_0-z_0^2$. 
\end{proof}

We will also  write (\ref{first}) in the form 
\begin{equation}\label{first2}
z(s)^2-a z(s)\cos\theta(s)+b\sin^2\theta(s)+az_0-z_0^2=0.
\end{equation}

We show that our solutions have symmetries at the critical points of the function $z$. 
Exactly, 

\begin{theorem}[Symmetry] Let $\alpha(s)=(x(s),z(s))$ be a profile curve of 
(\ref{eq1}) of a hyperbolic rotational LW-surface. If for some 
$s_1\in\r$, $\sin\theta(s_1)=0$, then $\alpha$ is symmetric with respect to the 
line $x=x(s_1)$.
\end{theorem}

\begin{proof} By hypothesis, $z'(s_1)=0$. We know that the functions 
$(x,z,\theta)$ are solutions of (\ref{eq1}). 
Without loss of generality, we assume 
that $x(s_1)=0$. Then it suffices to show
\begin{eqnarray*}
x(s_1-s)&=&-x(s+s_1)\\
z(s_1-s)&=&z(s+s_1)\\
\theta(s_1-s)&=&-\theta(s+s_1)
\end{eqnarray*}
But each pair of three functions is a solution of the same differential equations system for 
the same initial conditions. The uniqueness of solutions finishes the proof. 
\end{proof}

\begin{figure}[htbp]\begin{center}
\includegraphics[width=11cm]{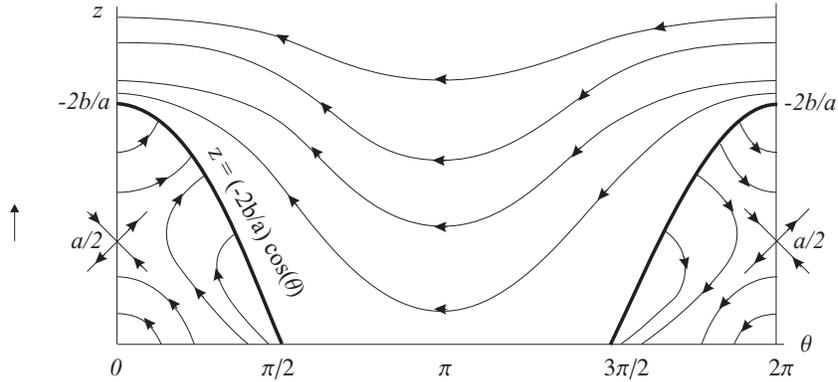}\end{center}
\caption{The phase portrait of the differential equations system (\ref{eq1}).}\label{fig10}
\end{figure}

We end this section describing the phase portrait of the differential equations system 
(\ref{eq1}), which allows a good understanding of the evolution of this 
system. Due to the periodicity of the functions cosine and sine, it suffices the 
study of the system (\ref{eq1}) for $\theta\in [0,2\pi]$. We project  the vector 
field (\ref{eq1}) into the $(\theta,z)$-plane, that is, 
\begin{eqnarray*}
\theta'(s)&=& \frac{a\cos\theta(s)-2z(s)}{a z(s)+2b \cos\theta(s)}\\
z'(s)&=& \sin\theta(s)
\end{eqnarray*}
In the region $[0,2\pi]\times\{(\theta,z);z>0\}$, the vector field has exactly two 
singularities at the points $(0,a/2)$ and $(2\pi,a/2)$ and there is a curve where 
the vector field $(\theta',z')$  is not defined, namely, $\{(\theta,z=(-2b/a)\cos\theta);0<\theta<\pi/2, 
3\pi/2< \theta <2\pi\}$. Both singularities are  saddle points, because 
the two eigenvalues of the linearization at the singularities have opposite signs.
 Exactly, the two eigenvalues are $\pm 2/\sqrt{-\Delta}$.  The singularities are of hyperbolic type. Moreover, 
 the eigenvectors of the  linearized system are  orthogonal. 
A sketch of the phase portrait of the system (\ref{eq1}) appears in Figure \ref{fig10}. 

Between the qualitative properties that we deduce from the phase portrait, we have the 
following ones depending on the range of the initial value $z_0$ in (\ref{eq2}):
\begin{enumerate}
\item If $0<z_0</2$, the generating curve $\alpha$ is defined in some bounded interval. Since the derivatives are bounded, it is necessary that the function $z$ vanishes at the extremum of the interval.
\item If $z_0=a/2$, we have a stationary solution, that is, $\alpha$ is a 
straight-line.
\item If $a/2<z_0<-2b/a$, the solution exploits when $z$  is nearing to $-2b/a \cos\theta$. This means that the maximal domain of the solution is some bounded interval.
\item If $z_0>-2b/a$,  it follows that the profile curve $\alpha$ is periodic.
\end{enumerate}

 According to this classification, due to when both 
numerator and denominator of $\theta'$ in (\ref{eq1}) vanish, we distinguish in the next sections  four cases depending 
on the value of the initial value $z_0$, namely, 
(i) $0<z_0<a/2$; (ii) $z_0=a/2$; (iii) $a/2<z_0<-2b/a$ and (iv) $z_0>-2b/a$. 
We point out that if $z_0=-2b/a$, there are not solutions of (\ref{eq1})-(\ref{eq2}).

\section{The case $z_0<a/2$}\label{s1}

In this section we consider the first case to discuss about the value $z_0$ in (\ref{eq2}): 
we assume, 
\begin{equation}\label{c1}
0<z_0<\frac{a}{2}.
\end{equation}
Let $(x,z,\theta)$ the solution of (\ref{eq1})-(\ref{eq2}). First, we study the qualitative 
properties fo the curve $\alpha(s)$, specially about its curvature $\theta$, and next we 
shall summarize the obtained results.

The function $z=z(s)$ satisfies 
 $z(0)=z_0$ and $z'(0)=0$. Because, $z_0<a/2<-2b/a$,  $\theta'(0)<0$. Then 
$z''(0)<0$. It follows that $z'$ and $z''$ are negative for 
$s>0$  near to $0$. 

\begin{claim} The function $z$ satisfies that $z'(s)<0$ for any $s>0$.
\end{claim}

 On the contrary, 
let us denote $s_1>s_0$  the first point where $z'(s_1)=0$. In particular, 
 $z''(s_1)\geq 0$. The function $\theta'$ is negative in the interval $[0,s_1]$. The proof 
is as follows. At $s=0$, $\theta'(0)<0$. Suppose that   $s=\bar{s}$, $0<\bar{s}\leq s_1$, 
 is the first zero of $\theta'$. Define the function $f(z_0):=z_0^2-az_0-b$. Using that $\Delta<0$, one can show that $f$ is always positive with a unique 
 minimum at $z_0=a/2$. Then (\ref{first}) implies that
$$\frac{a^2+4b}{a^2}z(\bar{s})^2=-f(z_0)<f(\frac{a}{2})=\frac{a^2+4b}{4}.$$
Then $z(\bar{s})^2>a^2/4>z_0^2$, which it is a contradiction because $z$ is decreasing in 
the interval $[0,s_1]$. As conclusion, $\theta'<0$ in $[0,s_1]$. 

We return with the function $z''$. As $z'(s_1)=0$, then $\cos\theta(s_1)=\pm 1$, but 
$z''(s_1)=\theta'(s_1)\cos\theta(s_1)\geq 0$ implies that $\cos\theta(s_1)=-1$. 
By using (\ref{first2}), $z(s_1)^2+az(s_1)+a z_0-z_0^2=0$, 
 which implies that $z(s_1)=-z_0$ or $z(s_1)=-a+z_0<0$: contradiction. This proves the Claim. 

Once proved the Claim, we show that $\theta'\not=0$ for all $s$. If the numerator is 
zero for some $s=s_2$ in the first zero, then   $\cos\theta(s_2)=2z(s_2)/a$. 
We use (\ref{first}) to deduce that 
$$\frac{a^2+4b}{a^2}z(s_2)^2=-f(z_0)<-f(\frac{a}{2})=\frac{a^2+4b}{4}.$$
Then $z(s_2)^2>a^2/4$. But as $z$ is strictly decreasing, $z(s_2)^2<z_0^2$ and we 
get a contradiction with the fact that $z_0<a/2$. Therefore  $\theta'<0$, that is,  
$\theta$  a decreasing and negative function for $s>0$.

We show that $z''<0$ in all the domain. Assume that for some $s>0$,  $z''(s)=0$. Since $\theta'<0$, it follows 
that $\cos\theta(s)=0$. As $a \cos\theta(s)-2z(s)>0$ for all $s$, we obtain a contradiction.

We study the maximal domain of definition of the given solution 
$(x,z,\theta)$. From the above reasoning, we have two possibilities: either $z=0$ at some point, that is, the curve $\alpha$ meets the $x$-axis,  or the maximal interval $[0,s_1)$ satisfies $s_1<\infty$ and $\lim_{s\rightarrow s_1}\theta'(s):=\theta_1=-\infty$ with 
$\lim_{s\rightarrow s_1}z(s_1):=z_1>0$. 
We see that the last one is impossible. From here, 
$a z_1+2b \cos\theta_1=0$. From (\ref{first}), 
$$(a^2+4b)z_1^2-4b f(z_0)=0\hspace*{.5cm}\mbox{or}\hspace*{.5cm}
z_1^2=\frac{4b f(z_0)}{a^2+4b}>\frac{4b f(a/2)}{a^2+4b}=-b.$$
Since $z(s)$ is a decreasing function on $s$, we deduce that
$$-b<z_1^2<z_0^2<\frac{a^2}{4},$$
which yields a contradiction. As conclusion, the function $z(s)$ 
vanishes at a first point $s_1$: $z(s_1)=0$. 
At this time, we know that the curve $\alpha$ can not be defined beyond $s=s_1$. 
Moreover,  $x'(s)\not =0$ for any $s$.

\begin{theorem} \label{t1} 
Let $\alpha=\alpha(s)=(x(s),0,z(s))$ be the profile curve of a hyperbolic rotational LW-surface $S$. Assume that the initial condition on $z_0$ satisfies (\ref{c1}).  Then 
\begin{enumerate}
\item The curve $\alpha$ is a graph on some bounded interval $(-x_1,x_1)$ of the $x$-axis, 
hence that $\alpha$ is embedded.
\item The curve 
$\alpha$  intersects the axis of rotation  at $x=\pm x_1$.
\item The curve $\alpha$ is concave, with exactly one maximum.
\end{enumerate}
\end{theorem}

\begin{theorem} \label{t12} Let $S$ be a hyperbolic rotational LW-surfaces whose profile 
curve $\alpha$ satisfies the hypothesis of the above Theorem. Then $S$ has the 
following properties:
\begin{enumerate}
\item The surface is embedded.
\item The Gaussian curvature of $S$ is positive. 
\item The surface $S$ can not to extend to be complete. Moreover, $S$ has  exactly two singular points which coincide with  the intersection of $S$ with its axis of rotation.
\end{enumerate}
\end{theorem}

 \begin{proof}  We only point out that $K$ is positive because 
 in the expression (\ref{hk}) for the Gauss curvature $K$, $\cos\theta>0$ and $\theta'<0$.
\end{proof} 

In  Figure \ref{fig1} (a), we show a picture of a curve $\alpha$ considered in this section. The surfaces that generate behave like the surfaces of revolution with positive constant Gauss 
curvature. See \cite{da}, \cite{ei}. 

\section{The case $z_0=a/2$:  cylinders}\label{s2}

\begin{theorem} \label{t2} 
Let $\alpha=\alpha(s)=(x(s),0,z(s))$ be the profile curve of a hyperbolic rotational LW-surface $S$. Assume that the initial condition on $z_0$ satisfies 
\begin{equation}\label{c2}
z_0=\frac{a}{2}.
\end{equation}
 Then $\alpha$ is a horizontal straight-line and $S$ is a right  cylinder.
\end{theorem}

\begin{proof} It is immediate that $x(s)=s$, $z(s)=a/2$ and $\theta(s)=0$ is 
the solution of (\ref{eq1})-(\ref{eq2}).
\end{proof}

\begin{figure}[htbp]
\includegraphics[width=6cm]{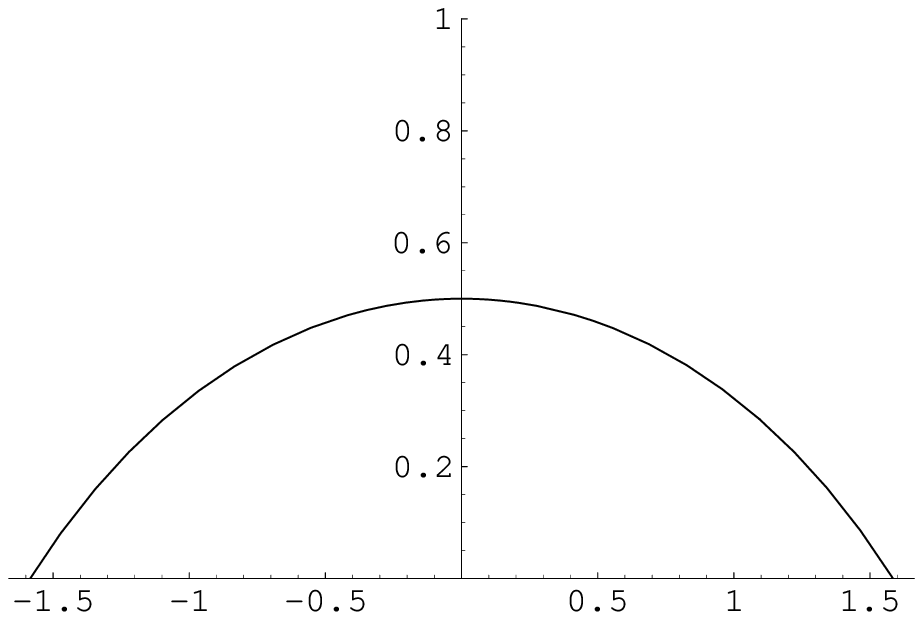}\hspace*{2cm}
\includegraphics[width=6cm]{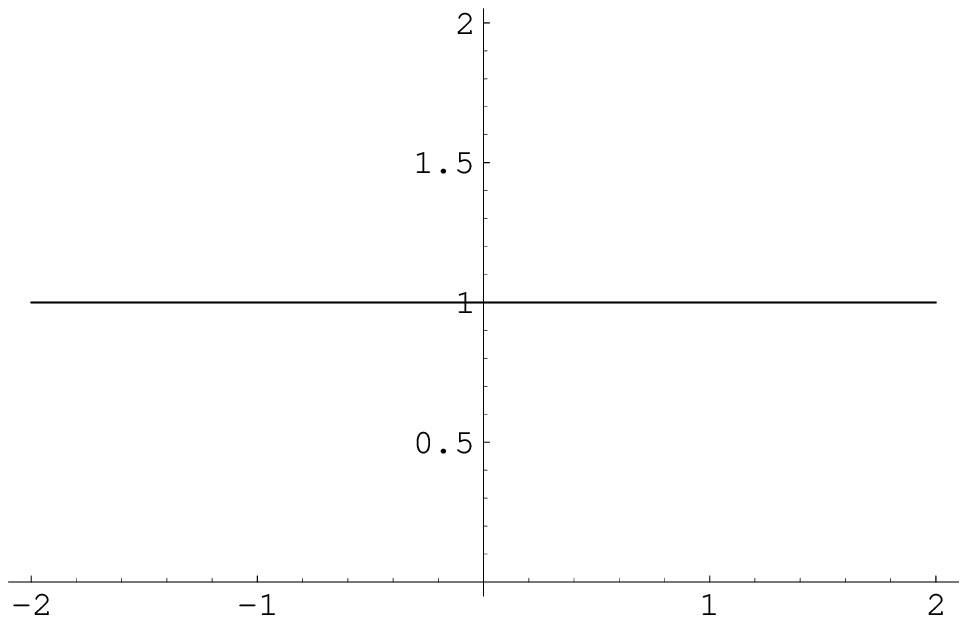}\\
\hspace*{2.7cm}(a)\hspace*{7.8cm}(b)
\caption{Two profile curves corresponding to rotational hyperbolic LW-surfaces. 
We assume that $a=-b=2$ in (\ref{w1}). (a) Case  $z_0=0.5$. The maximal domain of 
the solution is approximately $(-1.69,1.69)$. The curve is concave with one maximum; (b) Case $z_0=1$. The solution is a 
horizontal straight-line and the surface that generates is a right-cylinder.}\label{fig1}
\end{figure}

See Figure \ref{fig1} (b).

\section{The case $a/2<z_0<-2b/a$}\label{s3}

We study the properties of the solutions of (\ref{eq1})-(\ref{eq2}) when the initial 
condition on $z_0$  satisfies 
\begin{equation}\label{c3}
\frac{a}{2}<z_0<\frac{-2b}{a}.
\end{equation}
 Under this 
hypothesis, the value of $\theta'$ at $s=0$ in the expression (\ref{eq1}) is 
positive: exactly,  both numerator and denominator are negative. 

As $\theta'(0)>0$, 
the function $z$ is strictly increasing at $s=0$. 
We prove that $z'(s)>0$ for any $s>0$. On the contrary, if $s_1$ is the first point where 
$z'(s_1)=0$, we have $z''(s_1)\leq 0$ and $z$ is strictly increasing 
in $[0,s_1)$. The numerator of $\theta'$ in (\ref{eq1})
does not vanish in the interval $[0,s_1]$ because in such case, if  $a\cos\theta(\bar{s})-2z(\bar{s})=0$ for some $\bar{s}$, $0<\bar{s}\leq s_1$, then 
$z(\bar{s})\leq a/2$: contradiction, because $z_0<z(\bar{s})$. 
On the other hand, in the interval $[0,s_1]$, the function $\cos\theta$ 
does not vanish in $s\in [0,s_1]$, since if $\cos\theta(s)=0$ for some $s$, then $\theta'(s)=-a/2<0$. As conclusion,  $z''(s_1)=\theta'(s_1)\cos\theta(s_1)>0$: contradiction.

With the same reasoning, one shows  that the functions $z'$ and $z''$ are positive in its 
maximal domain $[0,s_1)$. We prove that $s_1$  must be finite. The proof is by contradiction.
Assume $s_1=\infty$. Then 
\begin{equation}\label{limit}
\lim_{s\rightarrow\infty}z(s)=\infty.
\end{equation}
At $s=0$, the value of the denominator  of $\theta'$  in (\ref{eq1}) is $az_0+2b$, 
which it is negative. However, using (\ref{limit}), $az(s)+2b\cos\theta(s)\rightarrow\infty$ as $s\rightarrow \infty$. This means that the denominator  of $\theta'$ 
in (\ref{eq1}) must   vanish  
at some point: contradiction.

\begin{theorem} \label{t3} 
Let $\alpha=\alpha(s)=(x(s),0,z(s))$ be the profile curve of a hyperbolic rotational LW-surface $S$. Assume that the initial condition on $z_0$ satisfies (\ref{c3}).  Then 
\begin{enumerate}
\item The curve $\alpha$ is a graph on some bounded interval $(-x_1,x_1)$ of the $x$-axis, 
hence that $\alpha$ is embedded.
\item The curve $\alpha$ is convex, with exactly one minimum.
\end{enumerate}
\end{theorem}

\begin{figure}[htbp]
\includegraphics[width=6cm]{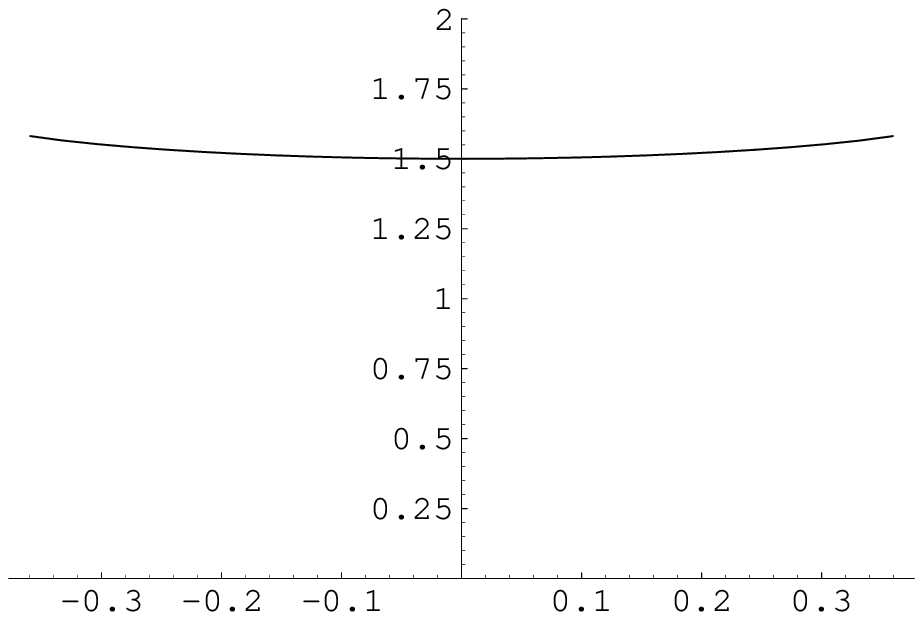}\hspace*{2cm}
\includegraphics[width=6cm]{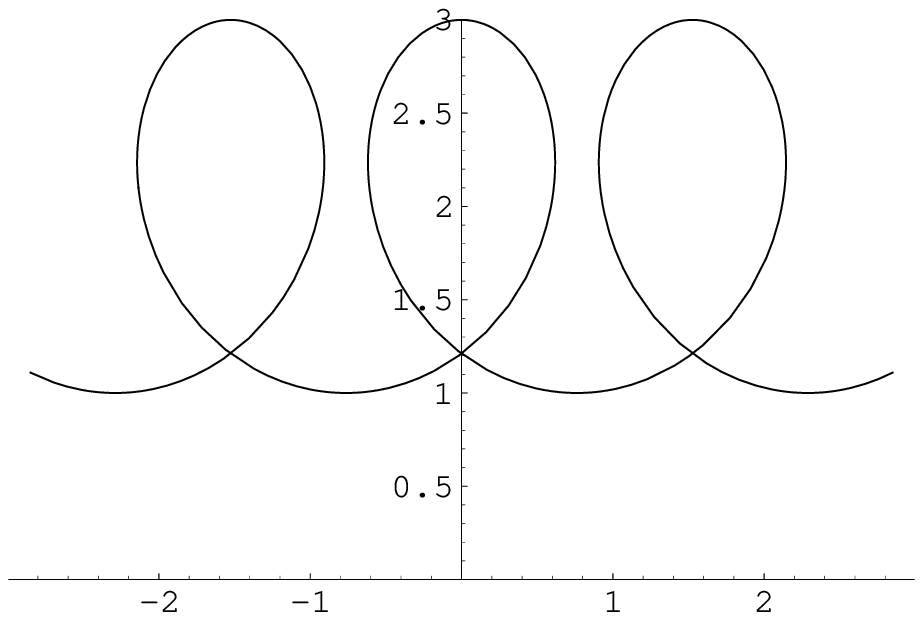}\\
\hspace*{2.7cm}(a)\hspace*{7.8cm}(b)
\caption{Two profile curves corresponding to rotational hyperbolic LW-surfaces. 
We assume that $a=-b=2$ in (\ref{w1}). (a) Case  $z_0=1.5$. The domain of 
the solution is approximately $(-0.372,0.372)$. Here $\alpha$ is convex with one minimum; 
(b) Case $z_0=3$. The curve $\alpha$ is periodic  with self-intersections. }\label{fig2}
\end{figure}

In  Figure \ref{fig2} (a), we present the profile curve $\alpha$ of a surface corresponding to the case studied in this section. As both $\cos\theta$ and $\theta'$ are 
positive functions, the Gaussian curvature $K$ in (\ref{hk}) is negative. 

\begin{theorem} \label{t32}
Let $S$ be a hyperbolic rotational LW-surfaces whose profile 
curve $\alpha$ satisfies the hypothesis of the above Theorem. Then $S$ has the 
following properties:
\begin{enumerate}
\item The surface is embedded.
\item The Gaussian curvature of $S$ is negative. 
\item The surface $S$  can not to extend to be  complete.
\end{enumerate}
\end{theorem}

\section{The case $z_0>-2b/a$: complete and periodic surfaces}\label{s4}
In this section, we study the initial value problem (\ref{eq1})-(\ref{eq2}) with 
the assumption 
\begin{equation}\label{c4}
z_0>\frac{-2b}{a}.
\end{equation}

From (\ref{first}), we write the function $z=z(s)$ as 
\begin{equation}\label{zeta}
z(s)=\frac12\left(a\cos\theta(s)+\sqrt{(a^2+4b)\cos^2\theta(s)+4(z_0^2-az_0-b)}\right).
\end{equation}

\begin{lemma}\label{le1}
 The maximal interval of the solution $(x,z,\theta)$ of  (\ref{eq1})-(\ref{eq2}) is $\r$.
\end{lemma}

\begin{proof} The result follows if we prove that the derivatives $x', z'$ and $\theta'$ are bounded. In view of (\ref{eq1}), it suffices to show it for $\theta'$:  we 
shall find negative numbers $m$ and $M$ such that $m\leq \theta'(s)\leq M<0$ for all $s$. We 
note that $\theta'(0)<0$.

First, we show the existence of constants $\delta_1$ and $\eta_1$, with 
 $\eta_1<0<\delta_1$, independent on $s$, such that 
\begin{equation}\label{delta}
az(s)+2b\cos\theta(s)\geq\delta_1\hspace*{.5cm}\mbox{and}\hspace*{.5cm} 
a\cos\theta(s)-2z(s)\geq\eta_1.
\end{equation}
Once proved this, it follows from (\ref{eq1}) that
\begin{equation}\label{m1}
\theta'(s)\geq \frac{\eta_1}{\delta_1}:=m.
\end{equation}
Because the function $f(z_0)$ is strictly increasing on $z_0$ for 
$z_0>a/2$,  there exists  $\epsilon>0$  such that 
$$z_0^2-a z_0-b=f(\frac{-2b}{a})+\epsilon=\frac{b(a^2+4b)}{a^2}+\epsilon.$$
From (\ref{zeta}), 
$$z\geq\frac12\left(a\cos\theta+\sqrt{(a^2+4b)\cos^2\theta+\frac{4b}{a^2}(a^2+4b)+4\epsilon}\right)\geq \frac12 (a\cos\theta-\frac{a^2+4b}{a}+\epsilon'),$$
for a certain positive number $\epsilon'$. By using the hyperbolicity condition 
$\Delta<0$, we have 
$$az(s)+2b\cos \theta(s)\geq \frac{a^2+4b}{2}(\cos\theta(s)-1)+\frac{a}{2}\epsilon'\geq \frac{a}{2}\epsilon':=\delta_1.$$
On the other hand, and using (\ref{zeta}) again
$$a\cos\theta(s)-2z(s)\geq-\sqrt{(a^2+4b)\cos^2(s)\theta+4f(z_0)}\geq -2\sqrt{f(z_0)}:=\eta_1.$$

We now obtain the upper bound for $\theta'$, $\theta'\leq M$. We prove that 
there exist $\delta_2,\eta_2$, with $\eta/2<0<\delta_2$ such that
\begin{equation}\label{delta2}
az(s)+2b\cos\theta(s)\leq\delta_2\hspace*{.5cm}\mbox{and}\hspace*{.5cm} 
a\cos\theta(s)-2z(s)\leq\eta_2.
\end{equation}
Using (\ref{zeta}), 
\begin{eqnarray*}
az(s)+2b\cos\theta(s)&=&\frac12\left((a^2+4b)\cos\theta(s)+a\sqrt{(a^2+4b)\cos^2\theta(s)+4f(z_0)}\right)\\
&\leq& a\sqrt{f(z_0)}:=\delta_2.
\end{eqnarray*}
On the other hand, 
\begin{eqnarray*}
a\cos\theta(s) -2z(s)&=&-\sqrt{(a^2+4b)\cos^2\theta(s)+4f(z_0)}\\
&\leq &-\sqrt{(a^2+4b)+4f(-2b/a)}:=\eta_2.
\end{eqnarray*}
Hence, we deduce from (\ref{eq1}) that 
\begin{equation}\label{ine}
\theta'(s)\leq \frac{\eta_2}{\delta_2}:=M.
\end{equation}
The inequalities (\ref{m1}) and (\ref{ine}) concludes the proof of the lemma.
\end{proof}

As a consequence of the proof of Lemma \ref{le1}, the graphic
 of the function $\theta$ lies between two slopped straight-lines. Since the derivative of $\theta$ is negative, we obtain that $\theta$ is strictly decreasing with 
$$\lim_{s\rightarrow\infty}\theta(s)=-\infty.$$
Put  $T>0$ the first number such that $\theta(T)=-2\pi$. We prove that $\alpha$ is a periodic curve.

\begin{lemma} \label{le2}
Under the hypothesis of this section and with the above notation, we have:
\begin{eqnarray*}
x(s+T)&=&x(s)+x(T)\\
z(s+T)&=&z(s)\\
\theta(s+T)&=&\theta(s)-2\pi
\end{eqnarray*}
\end{lemma}

\begin{proof}
This is a consequence of the uniqueness of solutions of (\ref{eq1})-(\ref{eq2}). We only have to show 
that $z(T)=z_0$. But this direct from (\ref{zeta}), the assumption (\ref{c4}) and 
that $a/2<-2b/a$ by the hyperbolicity condition $\Delta<0$.
\end{proof}

As conclusion of Lemmas \ref{le1} and \ref{le2},  we describe 
the behavior of the coordinates functions of the profile curve $\alpha$ under the 
assumption (\ref{c4}). See  Figure \ref{fig2} (b). 
Due to the monotonicity of 
$\theta$, let $T_1, T_2$ and $T_3$ be the points in the period $[0,T]$ such that the function 
$\theta$ takes the values $-\pi/2, -\pi$ and $-3\pi/2$ respectively. In view of the variation of the angle $\theta$ with the time coordinate $s$, it is easy to verify the following  Table:
\vspace*{4mm}

\begin{center}
\begin{tabular}{|c|c|c|c|}\hline
$s$ & $\theta$ & $x(s)$ & $z(s)$\\ \hline
$[0,T_1]$ & $[0,\frac{-\pi}{2}]$ & \mbox{ increasing} 
& \mbox{ decreasing}\\ \hline
$[T_1,T_2]$ & $[\frac{-\pi}{2},-\pi]$ &\mbox{ decreasing} &
\mbox{ decreasing} \\ \hline
$[T_2,T_3]$ & $[-\pi,\frac{-3\pi}{2}]$ &\mbox{ decreasing} &
\mbox{ increasing} \\ \hline
$[T_3,T]$ & $[\frac{-3\pi}{2},-2\pi]$ &\mbox{ increasing} & 
\mbox{ increasing}\\ 
\hline
\end{tabular}
\end{center}
\vspace*{4mm}

\begin{theorem} \label{t4} 
Let $\alpha=\alpha(s)=(x(s),0,z(s))$ be the profile curve of a hyperbolic rotational LW-surface $S$. Assume that the initial condition on $z_0$ satisfies (\ref{c4}).  Then 
\begin{enumerate}
\item The curve $\alpha$ is invariant by the group of translations in the $x$-direction given by the vector $(x(T),0,0)$. 
\item In each period of $z$, the curve $\alpha$ presents a maximum at $s=0$ 
and a minimum at $s=T_2$.  Moreover, $\alpha$ is symmetric with respect to the vertical line at $x=0$ and $x=x(T_2)$.
\item The height function of  $\alpha$, namely, $z=z(s)$,  is periodic. 
\item The curve $\alpha$ has self-intersections and its curvature has constant sign.
\item The part of $\alpha$ between the maximum and the minimum satisfies that the function $z(s)$ is 
strictly decreasing with exactly one vertical point. Between this
minimum and the next maximum, $z=z(s)$ is strictly increasing with exactly one vertical point.
\item The velocity $\alpha'$  turns around the origin.
\end{enumerate}
\end{theorem}

\begin{theorem} \label{t42} Let $S$ be a hyperbolic rotational LW-surfaces whose profile 
curve $\alpha$ satisfies the hypothesis of the above Theorem. Then $S$ has the 
following properties:
\begin{enumerate}
\item The surface has self-intersections.
\item The surface is periodic with infinite vertical symmetries.
\item The surface is  complete.
\item The part of $\alpha$ between two consecutive vertical points and containing a maximum corresponds with points 
of $S$ with positive Gaussian curvature; on the contrary, if this part contains a minimum,  the Gaussian curvature 
of the corresponding points of $S$ is negative on the surface.
\end{enumerate}
\end{theorem}

\begin{corollary}  Let $\alpha(s)=(x(s),0,z(s))$ be the profile curve of a hyperbolic rotational LW-surface. Assume that $\alpha$ is the solution of (\ref{eq1})-(\ref{eq2}) 
where $z_0>-2b/a$. Then the graphic of $\alpha$ lies between the lines $z=z_0-a$ and 
$z=z_0$, that is, 
$$z_0-a\leq z(s)\leq z_0,$$
where $z(s)$ reaches the minimum and the maximum values  in a discret set of points.
\end{corollary}

\begin{proof} Because $\theta\rightarrow-\infty$, the minimum and the maximum of the function 
$z(s)$ reach   at those 
points with $\cos\theta=-1$ and $\cos\theta=1$ respectively. 
The estimate is obtained from (\ref{zeta}).
\end{proof}

As it was announced in the Introduction and with the objective to distinghish with the surfaces of negative constant Gaussian curvature, 
 we stand out from Theorem \ref{t42} the following

\begin{corollary}\label{co}There exists a one-parameter family of 
rotational hyperbolic linear Weingarten 
surfaces that 
are complete and with self-intersections in $\r^3$. Moreover, these surfaces are periodic.
\end{corollary}




\begin{thebibliography}{99}

\bibitem{be} F. B. Brito, R. Sa Earp, 
On the structure of certain Weingarten surfaces with
boundary a circle, Ann. Fac. Sci. Toulouse, 6 (1997) 243--256,

\bibitem{ch1} S. S. Chern,
 Some new characterization of the Euclidean sphere, Duke Math. J., 12 (1945) 279--290

\bibitem{ch2}  S. S. Chern, 
On special $W$-surfaces, Proc. Amer. Math. Soc., 6 (1955) 783--786.


\bibitem{cft} A. V. Corro, W. Ferreira and K. Tenenblat,
Ribaucour transformations for constant mean curvature and linear Weingarten surfaces,
Pacific J. Math. 212 (2003) 265--296.

\bibitem{da} G. Darboux, 
Le\c{c}ons sur la th\'eorie des surfaces et les applications
geometrique du calcul infinitesimal, t. 1-4, Paris, Gauthier-Villars,
1877--1896.

\bibitem{ei} L. P.  Eisenhart,
 A treatise on the differential geometry of curves and
surfaces, Dover Publ., New York, 1909.


\bibitem{gmm}  J. A. G\'alvez, A.  Mart\'{\i}nez, F.  Milá\'an,
Linear Weingarten surfaces in $\r^3$,   Monatsh. Math. 138 (2003) 133--144.

\bibitem{hw} P. Hartman, W.  Winter,
 Umbilical points and W-surfaces, Am.  J. Math. 76 (1954) 502-–508.

\bibitem{ho}  H. Hopf,
 \"{U}ber Fl\"{a}chen mit einer Relation zwischen den Hauptkr\"{u}mmungen, 
Math. Nachr. 4 (1951) 232--249.

\bibitem{ho2} H. Hopf,
 Differential Geometry in the Large, Lecture Notes in Math, vol. 1000 (1983)
Berlin, Springer-Verlag.

\bibitem{ks}  W. K\"{u}hnel, M.  Steller, 
On closed Weingarten surfaces,  Monatsh. Math. 146 (2005) 113--126.


\bibitem{rs} H. Rosenberg,   R. Sa Earp,
The geometry of properly embedded special surfaces in $\r^3$; e.g., surfaces satisfying $aH+bK=1$, where $a$ and $b$ are positive, 
 Duke Math. J. 73 (1994)  291-–306.


\bibitem{vo}  K. Voss,
\"{U}ber geschlossene Weingartensche Fl\"{a}chen,
Math. Annalen, 138 (1959) 42--54.


\bibitem{we1} J. Weingarten,
Ueber eine Klasse auf einander abwickelbarer Fl\"{a}achen,
J. Reine Angew. Math. 59 (1861)  382--393.

 \bibitem{we2} J. Weingarten,  
Ueber die Fl\"{a}chen, derer Normalen eine gegebene Fl\"{a}che ber\"{u}hren, 
 J. Reine Angew. Math.  62 (1863) 61-–63.

\end{thebibliography}
\end{document}